\documentclass{amsart}
\usepackage[enableskew,vcentermath]{youngtab}

\newtheorem{theorem}{Theorem}
\newtheorem{lemma}[theorem]{Lemma}

\newtheorem{corollary}[theorem]{Corollary}
\newtheorem{remark}[theorem]{Remark}

\newcommand{\twop}{2'}
\newcommand{\threep}{3'}
\newcommand{\fourp}{4'}

\newcommand{\wh}{\widehat{w}}
\newcommand{\ear}{\overline{e}}
\DeclareMathOperator{\maj}{maj}
\newcommand{\shsyt}{SYT'(2n-1, 2n-3, \ldots, 1)}
\newcommand{\syt}{SYT(n^n)}

\begin{document}

\author[T. K. Petersen]{T. Kyle Petersen}

\author[L. Serrano]{Luis Serrano}

\title[Cyclic sieving in the hyperoctahedral group]{Cyclic sieving for longest reduced words \\ in the hyperoctahedral group}

\begin{abstract}
We show that the set $R(w_0)$ of reduced expressions for the longest element in the hyperoctahedral group exhibits the cyclic sieving phenomenon. More specifically, $R(w_0)$ possesses a natural cyclic action given by moving the first letter of a word to the end, and we show that the orbit structure of this action is encoded by the generating function for the major index on~$R(w_0)$.
\end{abstract}

\maketitle

\section{Introduction and main result}

Suppose we are given a finite set $X$, a finite cyclic group $C = \langle \omega \rangle$ acting on $X$, and a polynomial $X(q) \in \mathbb{Z}[q]$ with integer coefficients.  Following Reiner, Stanton, and White \cite{RSW}, we say that the triple $(X, C, X(q))$ exhibits the \emph{cyclic sieving phenomenon} (CSP) if for every integer $d \geq 0$, we have that $|X^{\omega^d}| = X(\zeta^d)$ where $\zeta \in \mathbb{C}$ is a root of unity of multiplicitive order $|C|$ and $X^{\omega^d}$ is the fixed point set of the action of the power $\omega^d$. The sizes of the fixed point sets determine the cycle structure of the canonical image of $\omega$ in the group of permutations of $X$, $S_X$.  Therefore, to find the cycle structure of the image of any bijection $\omega: X \rightarrow X$, it is enough to determine the order of the action of $\omega$ on $X$ and find a polynomial $X(q)$ such that $(X, \langle \omega \rangle, X(q))$ exhibits the CSP. 

The cyclic sieving phenomenon has been demonstrated in a variety of contexts. The paper of Reiner, Stanton, and White \cite{RSW} itself includes examples involving noncrossing partitions, triangulations of polygons, and cosets of parabolic subgroups of Coxeter groups. An example of the CSP with standard Young tableaux is due to Rhoades \cite{Rh} and will discussed further in Section \ref{sec:Rhoades}. Now we turn to the CSP of interest to this note.

Let $w_0 = w_0^{(B_n)}$ denote the longest element in the type $B_n$ Coxeter group. Given generating set $S = \{s_1,\ldots,s_n\}$ for $B_n$, ($s_1$ being the ``special" reflection), we will write a reduced expression for $w_0$ as a word in the subscripts. For example, $w_0^{(B_3)}$ can be written as \[ s_1s_2s_1s_3s_2s_3s_1s_2s_3;\] we will abbreviate this product by $121323123$. It turns out that if we cyclically permute these letters, we always get another reduced expression for $w_0$. Said another way, $s_i w_0 s_i = w_0$ for $i = 1,\ldots,n$. The same is not true for longest elements of other classical types. In type A, we have $s_i w_0^{(A_n)} s_{n+1-i} = w_0^{(A_n)}$, and for type D, 
\[ w_0^{(D_n)} = \begin{cases} s_i w_0^{(D_n)} s_i & \mbox{if } n \mbox{ even or } i > 2,\\
s_i w_0^{(D_n)} s_{3-i} & \mbox{if } n \mbox{ odd and }  i=1,2.
\end{cases}\]

Let $R(w_0)$ denote the set of reduced expressions for $w_0$ in type $B_n$ and let $c: R(w_0) \to R(w_0)$ denote the action of placing the first letter of a word at the end. Then the orbit in $R(w_0^{(B_3)})$ of the word above is:
\[ \{121323123 \to 213231231 \to 132312312 \to 323123121 \to 231231213 \] \[\to 312312132 \to 123121323 \to 231213231 \to 312132312 \}. \]

As the length of $w_0$ is $n^2$, we clearly have $c^{n^2} = 1$, and the size of any orbit divides $n^2$. For an example of a smaller orbit, notice that the word $213213213$ has cyclic order 3.

For any word $w = w_1 \ldots w_l$, (e.g., a reduced expression for $w_0$), a \emph{descent} of $w$ is defined to be a position $i$ in which $w_i > w_{i+1}$.The \emph{major index} of $w$, $\maj(w)$, is defined as the sum of the descent positions. For example, the word $w = 121323123$ has descents in positions $2,4$, and $6$, so its major index is $\maj(w) = 2+4+6 = 12$. Let $f_n(q)$ denote the generating function for this statistic on words in $R(w_0)$:
\[
f_n(q) = \sum_{w \in R(w_0)} q^{\maj(w)}.
\]

The following is our main result.

\begin{theorem}\label{theorem:main}
The triple $(R(w_0),\langle c \rangle,X(q))$ exhibits the cyclic sieving phenomenon, where \[X(q) =q^{-n\binom{n}{2}}f_n(q).\]
\end{theorem}

For example, let us consider the case $n=3$. We have \begin{align*}
  X(q) &= q^{-9}\sum_{w \in R(w_0^{(B_3)})} q^{\maj(w)} \\
  &= 1+ q^2+2q^3+ 2q^4 + 2q^5+4q^6 + 3q^7+4q^8 +4q^9\\
 & \quad + 4q^{10}+3q^{11}+4q^{12} + 2q^{13} +2q^{14}+ 2q^{15}+ q^{16}+q^{18}.
 \end{align*}
Let $\zeta = e^{\frac{2\pi i}{9}}$. Then we compute:
\[ \begin{array}{ l c r }
X(1) = 42 & X(\zeta^3) = 6 & X(\zeta^6) = 6 \\
X(\zeta) = 0 & X(\zeta^4) = 0 & X(\zeta^7) = 0 \\
X(\zeta^2) = 0 & X(\zeta^5) = 0 & X(\zeta^8) = 0
\end{array} \]
Thus, the $42$ reduced expressions for $w_0^{(B_3)}$ split into two orbits of size three (the orbits of $123123123$ and $132132132$) and four orbits of size nine.

We prove Theorem \ref{theorem:main} by relating it to another instance of the CSP, namely Rhoades' recent (and deep) result \cite[Thm 3.9]{Rh} for the set $SYT(n^m)$ of rectangular standard Young tableaux with respect to the action of \emph{promotion} (defined in Section \ref{sec:prom}). To make the connection, we rely on a pair of remarkable bijections due to Haiman \cite{H1,H2}. The composition of Haiman's bijections maps to $R(w_0)$ from the set of square tableaux, $\syt$. In this note our main goal is to show that Haiman's bijections carry the orbit structure of promotion on $\syt$ to the orbit structure of $c$ on $R(w_0)$.

We conclude this section by remarking that this approach was first outlined by Rhoades \cite[Thm 8.1]{Rh}. One purpose of this article is to fill some nontrivial gaps in his argument. A second is to justify the new observation that the polynomial $X(q)$ can be expressed as the generating function for the major index on $R(w_0)$. We thank Brendon Rhoades for encouraging us to write this note. Thanks also to Kevin Dilks, John Stembridge, and Alex Yong for fruitful discussions on this and related topics, and to Sergey Fomin for comments on the manuscript.

\section{Promotion on standard Young tableaux}\label{sec:prom}

For $\lambda$ a partition, let $SYT(\lambda)$ denote the set of standard Young tableaux of shape $\lambda$. If $\lambda$ is a \emph{strict partition}, i.e., with no equal parts, then let $SYT'(\lambda)$ denote the set of standard Young tableaux of shifted shape $\lambda$. We now describe the action of \emph{jeu de taquin promotion}, first defined by Sch\"utzenberger \cite{Sch}.

We will consider promotion as a permutation of tableaux of a fixed shape (resp. shifted shape), $p: SYT(\lambda) \to SYT(\lambda)$ (resp. $p: SYT'(\lambda) \to SYT'(\lambda)$). Given a tableau $T$ with $\lambda \vdash n$, we form $p(T)$ with the following algorithm. (We denote the entry in row $a$, column $b$ of a tableau $T$, by $T_{a,b}$.)
\begin{enumerate}

\item Remove the entry 1 in the upper left corner and decrease every other entry by 1. The empty box is initialized in position $(a,b) = (1,1)$.

\item Perform jeu de taquin:
\begin{enumerate}

\item If there is no box to the right of the empty box and no box below the empty box, then go to 3).

\item If there is a box to the right or below the empty box, then swap the empty box with the box containing the smaller entry, i.e., $p(T)_{a,b} := \min\{T_{a,b+1}, T_{a+1,b}\}$. Set $(a,b) := (a',b')$, where $(a',b')$ are the coordinates of box swapped, and go to 2a).

\end{enumerate}

\item Fill the empty box with $n$.

\end{enumerate}
Here is an example:
\[T = \young(1248,367,5) \quad \mapsto  \quad \young(1367,258,4) = p(T).
\]

As a permutation, promotion naturally splits $SYT(\lambda)$ into disjoint orbits. For a general shape $\lambda$ there seems to be no obvious pattern to the sizes of the orbits. However, for certain shapes, notably Haiman's ``generalized staircases" more can be said \cite{H2} (see also Edelman and Greene \cite[Cor. 7.23]{EG}). In particular, rectangles fall into this category, with the following result.

\begin{theorem}[\cite{H2}, Theorem 4.4] \label{thm:ordern}
If $\lambda \vdash N = bn$ is a rectangle, then $p^N(T) = T$ for all $T \in SYT(\lambda)$.
\end{theorem}

Thus for $n \times n$ square shapes $\lambda$, $p^{n^2} =1$ and the size of every orbit divides $n^2$. With $n=3$, here is an orbit of size 3: 
\begin{equation}\label{eq:3orb} 
\young(125,368,479) \to \young(147,258,369) \to \young(136,247,589)\to \cdots.
\end{equation}

There are $42$ standard Young tableaux of shape $(3,3,3)$, and there are $42$ reduced expressions in the set $R(w_0^{(B_3)})$. Stanley first conjectured that $R(w_0)$ and $SYT(n^n)$ are equinumerous, and Proctor suggested that rather than $SYT(n^n)$, a more direct correspondence might be given with $SYT'(2n-1,2n-3,\ldots,1)$, that is, with shifted standard tableaux of ``doubled staircase" shape. (That the squares and doubled staircases are equinumerous follows easily from hook length formulas.)

Haiman answers Proctor's conjecture in such a way that the structure of promotion on doubled staircases corresponds precisely to cyclic permutation of words in $R(w_0)$ \cite[Theorem 5.12]{H2}. Moreover, in \cite[Proposition 8.11]{H1}, he gives a bijection between standard Young tableaux of square shape and those of doubled staircase shape that (as we will show) commutes with promotion.

As an example, his bijection carries the orbit in \eqref{eq:3orb} to this shifted orbit: \[ \young( 12458,:369,::7) \to \young( 12347,:568,::9) \to \young(12369,:457,::8) \to \cdots.\] Both of these orbits of tableaux correspond to the orbit of the reduced word $132132132$.

\section{Haiman's bijections}

We first describe the bijection between reduced expressions and shifted standard tableaux of doubled staircase shape. This bijection is described in Section 5 of~\cite{H2}.

Let $T$ in $SYT'(2n-1,2n-3,\ldots,1)$. Notice the largest entry in $T$, (i.e., $n^2$), occupies one of the outer corners. Let $r(T)$ denote the row containing this largest entry, numbering the rows from the bottom up. The \emph{promotion sequence} of $T$ is defined to be $\Phi(T)=r_1\cdots r_{n^2}$, where $r_i = r(p^i(T))$. Using the example above of \[T = \young( 12458,:369,::7),\] we see $r(T) = 2, r(p(T)) = 1$, $r(p^2(T)) = 3$, and since $p^3(T) = T$, we have \[ \Phi(T) = 132132132.\]

Haiman's result is the following.

\begin{theorem}[\cite{H2}, Theorem 5.12]\label{theorem:phiandccommute}
The map $T \mapsto \Phi(T)$ is a bijection $\shsyt \to R(w_0)$.
\end{theorem}

By construction, then, we have \[ \Phi(p(T)) = c(\Phi(T)),\] i.e., $\Phi$ is an orbit-preserving bijection \[ (\shsyt,p) \longleftrightarrow (R(w_0),c).\]

Next, we will describe the bijection \[H: \syt \to \shsyt\] between squares and doubled staircases. Though not obvious from the definition below, we will demonstrate that $H$ commutes with promotion. 

We assume the reader is familiar with the \emph{Robinson-Schensted-Knuth insertion} algorithm (RSK). (See \cite[Section 7.11]{Sta}, for example.) This is a map between words $w$ and pairs of tableaux $(P,Q) = (P(w), Q(w))$. We say $P$ is the \emph{insertion tableau} and $Q$ is the \emph{recording tableau}.

There is a similar correspondence between words $w$ and pairs of shifted tableaux $(P',Q') = (P'(w),Q'(w))$ called \emph{shifted mixed insertion} due to Haiman \cite{H1}. (See also Sagan \cite{Sa} and Worley \cite{W}.) Serrano defined a semistandard generalization of shifted mixed insertion in \cite{Ser}. Throughout this paper we refer to semistandard shifted mixed insertion simply as \emph{mixed insertion}. Details can be found in \cite[Section 1.1]{Ser}.

\begin{theorem}[\cite{Ser} Theorem 2.26]\label{theorem: bijection}
Let $w$ be a word. If we view $Q(w)$ as a skew shifted standard Young tableau and apply jeu de taquin to obtain a standard shifted Young tableau, the result is $Q'(w)$ (independent of any choices in applying jeu de taquin).
\end{theorem}

For example, if $w =332132121$, then
\[ (P,Q) = \left(\young(111,222,333),\young(125,368,479)\right),\] \[(P',Q') = \left(\young(111\twop\threep,:22\threep,::3), \young( 12458,:369,::7)\right).\] Performing jeu de taquin we see:
\[ \young(\hfil\hfil125,:\hfil368,::479) \to \young(\hfil\hfil125,:3468,::79) \to \young(\hfil1258,:346,::79) \to \young(12458,:369,::7).\] 

Haiman's bijection is precisely $H(Q) = Q'$. That is, given a standard square tableau $Q$, we embed it in a shifted shape and apply jeu de taquin to create a standard shifted tableau. That this is indeed a bijection follows from Theorem \ref{theorem: bijection}, but is originally found in \cite[Proposition 8.11]{H1}.

\begin{remark}
Haiman's bijection $H$ applies more generally between rectangles and ``shifted trapezoids", i.e., for $m\leq n$, we have $H: SYT(n^m) \to SYT'(n+m-1,n+m-3,\ldots, n-m+1)$. All the results presented here extend to this generality, with similar proofs. We restict to squares and doubled staircases for clarity of exposition.
\end{remark}

We will now fix the tableaux $P$ and $P'$ to ensure that the insertion word $w$ has particularly nice properties. We will use the following lemma.

\begin{lemma}[\cite{Ser}, Proposition 1.8]\label{lem:insert}
Fix a word $w$. Let $P = P(w)$ be the RSK insertion tableau and let $P'=P'(w)$ be the mixed insertion tableau. Then the set of words that mixed insert into $P'$ is contained in the set of words that RSK insert into $P$.
\end{lemma}

Now we apply Lemma \ref{lem:insert} to the word \[ w = \underbrace{n\cdots n}_n \cdots \underbrace{2\cdots 2}_n \underbrace{1\cdots 1}_n.\]

If we use RSK insertion, we find $P$ is an $n \times n$ square tableau with all 1s in row first row, all 2s in the second row, and so on. With such a choice of $P$ it is not difficult to show that any other word $u$ inserting to $P$ has the property that in any initial subword $u_1\cdots u_i$, there are at least as many letters $(j+1)$ as letters $j$. Such words are sometimes called (reverse) \emph{lattice words} or (reverse) \emph{Yamanouchi words}. Notice also that any such $u$ has $n$ copies of each letter $i$, $i=1,\ldots,n$. We call the words inserting to this choice of $P$ \emph{square words}.

On the other hand, if we use mixed insertion on $w$, we find $P'$ as follows (with $n=4$):\[ \young(1111\twop\threep\fourp,:222\threep\fourp,::33\fourp,:::4).\] In general, on the ``shifted half" of the tableau we see all 1s in the first row, all 2s in the second row, and so on. In the ``straight half" we see only prime numbers, with $2'$ on the first diagonal, $3'$ on the second diagonal, and so on. Lemma \ref{lem:insert} tells us that every $u$ that mixed inserts to $P'$ is a square word. But since the sets of recording tableaux for $P$ and for $P'$ are equinumerous, we see that the set of words mixed inserting to $P'$ is precisely the set of all square words.

\begin{remark}
Yamanouchi words give a bijection with square standard Young tableaux that circumvents insertion completely. In reading the word from left to right, if $w_i = j$, we put letter $i$ in the leftmost unoccupied position of row $n+1-j$. (See \cite[Proposition 7.10.3(d)]{Sta}.)
\end{remark}

We will soon characterize promotion in terms of operators on insertion words. First, some lemmas.

For a tableau $T$ (shifted or not) let $\Delta T$ denote the result of all but step (3) of promotion. That is, we delete the smallest entry and perform jeu de taquin, but we do not fill in the empty box. The following lemma says that, in both the shifted and unshifted cases, this can be expressed very simply in terms of our insertion word. The first part of the lemma is a direct application of the theory of jeu de taquin (see, e.g., \cite[A1.2]{Sta}); the second part is \cite[Lemma 3.9]{Ser}.

\begin{lemma}\label{lemma: recording}
For a word $w= w_1w_2\cdots w_l$, let $\wh = w_2 \cdots w_l$. Then we have
\[ Q(\wh) = \Delta Q(w),\] and
\[ Q'(\wh) = \Delta Q'(w).\]
\end{lemma}

The operator $e_j$ acting on words $w = w_1 \cdots w_l$ is defined in the following way. Consider the subword of $w$ formed only by the letters $j$ and $j+1$. Consider every $j+1$ as an opening bracket and every $j$ as a closing bracket, and pair them up accordingly. The remaining word is of the form $j^r(j+1)^s$. The operator $e_j$ leaves all of $w$ invariant, except for this subword, which it changes to $j^{r-1}(j+1)^{s+1}$. (This operator is widely used in the theory of \emph{crystal graphs}.)

As an example, we calculate $e_2(w)$ for the word $w = 3121221332$. The subword formed from the letters $3$ and $2$ is\[3\cdot2\cdot22\cdot332,\] which corresponds to the bracket sequence $()))(()$. Removing paired brackets, one obtains $))($, corresponding to the subword \[\cdot\cdot\cdot\cdot 22\cdot3\cdot\cdot .\] We change the last $2$ to a $3$ and keep the rest of the word unchanged, obtaining $e_2(w) = 31212\mathbf{3}1332$.

The following lemma shows that this operator leaves the recording tableau unchanged. The unshifted case is found in work of Lascoux, Leclerc, and Thibon \cite[Theorem 5.5.1]{LLT}; the shifted case follows from the unshifted case, and the fact that the mixed recording tableau of a word is uniquely determined by its RSK recording tableau (Theorem \ref{theorem: bijection}).

\begin{lemma}\label{lemma: operator}
Recording tableaux are invariant under the operators $e_i$. That is,
\[
Q(e_i(w)) = Q(w),
\]
and
\[ Q'(e_i(w)) = Q'(w).\]
\end{lemma}

Let $\ear = e_1 \cdots e_{n-1}$ denote the composite operator given by applying first $e_{n-1}$, then $e_{n-2}$ and so on. It is clear that if $w=w_1\cdots w_{n^2}$ is a square word, then $\ear(\wh)1$ is again a square word.

\begin{theorem}\label{theorem: commute}
Let $w=w_1\cdots w_{n^2}$ be a square word. Then,
\[ p(Q(w)) = Q(\ear(\wh)1),\]
and
\[ p(Q'(w)) = Q'(\ear(\wh)1).\] In other words,
 Haiman's bijection commutes with promotion: \[ p(H(Q)) = H(p(Q)).\]
\end{theorem}

\begin{proof}
By Lemma \ref{lemma: recording}, we see that $Q(\wh)$ is only one box away from $p(Q(w))$. Further, repeated application of Lemma \ref{lemma: operator} shows that \[Q(\wh) = Q(e_{n-1}(\wh)) = Q(e_{n-2}(e_{n-1}(\wh))) = \cdots = Q(\ear(\wh)).\] The same lemmas apply show $Q'(\ear(\wh))$ is one box away from $p(Q'(w))$. 

All that remains is to check that the box added by inserting 1 into $P(\ear(\wh))$ (resp. $P'(\ear(\wh))$) is in the correct position. But this follows from the observation that $\ear(\wh)1$ is a square word, and square words insert (resp. mixed insert) to squares (resp. doubled staircases).
\end{proof}

\section{Rhoades' result}\label{sec:Rhoades}

Rhoades \cite{Rh} proved an instance of the CSP related to the action of promotion on rectangular tableaux. His result is quite deep, employing Kahzdan-Lusztig cellular representation theory in its proof. 

Recall that for any partition $\lambda \vdash n$, we have that the standard tableaux of shape $\lambda$ are enumerated by the Frame-Robinson-Thrall \emph{hook length formula}:
\[
f^{\lambda} = |SYT(\lambda)| = \frac{n!}{\prod_{(i,j) \in \lambda} h_{ij}},
\]
where the product is over the boxes $(i,j)$ in $\lambda$ and $h_{ij}$ is the \emph{hook length} at the box $(i,j)$, i.e., the number of boxes directly east or south of the box $(i,j)$ in $\lambda$, counting itself exactly once.  To obtain the polynomial used for cyclic sieving, we replace the hook length formula with a natural $q$-analogue. First, recall that for any $n \in \mathbb{N}$, $[n]_q := 1 + q + \cdots + q^{n-1}$ and $[n]_q! := [n]_q [n-1]_q \cdots [1]_q$.

\begin{theorem}[\cite{Rh}, Theorem 3.9] \label{thm:cs}
Let $\lambda \vdash N$ be a rectangular shape and let $X = SYT(\lambda)$.  Let $C := \mathbb{Z} / N \mathbb{Z}$ act on $X$ via promotion.  Then, the triple $(X, C, X(q))$ exhibits the cyclic sieving phenomenon, where
\begin{equation*}
X(q) = \frac{[N]_q!}{\Pi_{(i,j) \in \lambda} [h_{ij}]_q}
\end{equation*}
is the $q$-analogue of the hook length formula.
\end{theorem}

Now thanks to Theorem \ref{theorem: commute} we know that $H$ preserves orbits of promotion, and as a consequence we see the CSP for doubled staircases.

\begin{corollary}\label{cor: sieving1}
Let $X = SYT'(2n-1,2n-3,\ldots,1)$, and let $C := \mathbb{Z}/ n^2\mathbb{Z}$ act on $X$ via promotion. Then the triple $(X, C, X(q))$ exhibits the cyclic sieving phenomenon, where
\[
X(q) = \frac{[n^2]_q!}{[n]_q^n\prod_{i=1}^{n-1} ([i]_q\cdot [2n-i]_q)^i}
\] is the $q$-analogue of the hook length formula for an $n\times n$ square Young diagram.
\end{corollary}

Because of Theorem \ref{theorem:phiandccommute} the set $R(w_0)$ also exhibits the CSP.

\begin{corollary}[\cite{Rh}, Theorem 8.1]\label{cor: sieving2}
Let $X = R(w_0)$ and let $X(q)$ as in Corollary \ref{cor: sieving1}.  Let $C := \mathbb{Z}/ n^2\mathbb{Z}$ act on $X$ by cyclic rotation of words. Then the triple $(X, C, X(q))$ exhibits the cyclic sieving phenomenon.
\end{corollary}

Corollary \ref{cor: sieving2} is the CSP for $R(w_0)$ as stated by Rhoades. This is nearly our main result (Theorem \ref{theorem:main}), but for the definition of $X(q)$.

In spirit, if $(X, C, X(q))$ exhibits the CSP, the polynomial $X(q)$ should be some $q$-enumerator for the set $X$. That is, it should be expressible as \[ X(q) = \sum_{x \in X} q^{s(x)},\] where $s$ is an intrinsically defined statistic for the elements of $X$. Indeed, nearly all known instances of the cyclic sieving phenomenon have this property.
For example, it is known (\cite[Cor 7.21.5]{Sta}) that the $q$-analogue of the hook-length formula can be expressed as follows: 
\begin{equation}\label{eq:hook}
f^{\lambda}(q) = q^{-\kappa(\lambda)} \sum_{T \in SYT(\lambda)} q^{\maj(T)},
\end{equation} where $\kappa(\lambda_1,\ldots,\lambda_l) = \sum_{1\leq i \leq l} (i-1)\lambda_i$ and for a tableau $T$, $\maj(T)$ is the sum of all $i$ such that $i$ appears in a row above $i+1$. Thus $X(q)$ in Theorem \ref{thm:cs} can be described in terms a statistic on Young tableaux.

With this point of view, Corollaries \ref{cor: sieving1} and \ref{cor: sieving2} are aesthetically unsatisfying. Section \ref{sec:comb} is given to showing that $X(q)$ can be defined as the generating function for the major index on words in $R(w_0)$. It would be interesting to find a combinatorial description for $X(q)$ in terms of a statistic on $\shsyt$ as well, though we have no such description at present.

\section{Combinatorial description of $X(q)$}\label{sec:comb}

As stated in the introduction, we will show that \[X(q) = q^{-n\binom{n}{2}}\sum_{w \in R(w_0)} q^{\maj(w)}.\] If we specialize \eqref{eq:hook} to square shapes, we see that $\kappa(n^n) = n\binom{n}{2}$ and 
\[ X(q) = q^{-n\binom{n}{2}}\sum_{T \in \syt} q^{\maj(T)}.\] Thus it suffices to exhibit a bijection between square tableaux and words in $R(w_0)$ that preserves major index. In fact, the composition $\Psi := \Phi H$ has a stronger feature.

Define the \emph{cyclic descent set} of a word $w=w_1\cdots w_l$ to be the set \[ D(w) = \{ i : w_i > w_{i+1}\} \pmod l\] That is, we have descents in the usual way, but also a descent in position 0 if $w_l > w_1$. Then $\maj(w) = \sum_{i \in D(w)} i$. For example with $w=132132132$, $D(w) = \{0,2,3,5,6,8\}$ and $\maj(w) = 0 + 2 + 3+ 5+ 6+8 = 24$. 

Similarly, we follow \cite{Rh} in defining the cyclic descent set of a square (in general, rectangular) Young tableau. For $T$ in $\syt$, define $D(T)$ to be the set of all $i$ such that $i$ appears in a row above $i+1$, along with 0 if $n^2-1$ is above $n^2$ in $p(T)$. Major index is $\maj(T) = \sum_{i \in D(T)}i$. We will see that $\Psi$ preserves cyclic descent sets, and hence, major index. Using our earlier example of $w=132132132$, one can check that \[T = \Psi^{-1}(w) = \young(125,368,479)\] has $D(T)=D(w)$, and so $\maj(T) = \maj(w)$. 

\begin{lemma}\label{lem:descents}
Let $T \in \syt$, and let $w = \Psi(T)$ in $R(w_0)$. Then $D(T) = D(w)$.
\end{lemma}

\begin{proof}
First, we observe that both types of descent sets shift cyclically under their respective actions: \[ D(p(T)) = \{ i-1 \pmod{n^2} : i \in D(T)\}, \] and \[ D(c(w)) = \{ i-1 \pmod{n^2} : i \in D(w)\}.\] For words under cyclic rotation, this is obvious. For tableaux under promotion, this is a lemma of Rhoades \cite[Lemma 3.3]{Rh}. 

Because of this cyclic shifting, we see that $i \in D(T)$ if and only if $0 \in D(p^i(T))$. Thus, it suffices to show that $0 \in D(T)$ if and only if $0 \in D(w)$. (Actually, it is easier to determine if $n^2-1$ is a descent.)

Let $S = \Phi^{-1}(w)$ be the shifted doubled staircase tableau corresponding to $w$. We have $n^2-1 \in D(w)$ if and only if $n^2$ is in a higher row in $p^{-1}(S)$ than in $S$. But since $n^2$ occupies the same place in $p^{-1}(S)$ as $n^2-1$ occupies in $S$, this is to say $n^2-1$ is above $n^2$ in $S$. On the other hand, $n^2-1 \in D(T)$ if and only if $n^2-1$ is above $n^2$ in $T$. It is straightforward to check that since $S$ is obtained from $T$ by jeu de taquin into the upper corner, the relative heights of $n^2$ and $n^2-1$ (i.e., whether $n^2$ is below or not) are the same in $S$ as in $T$. This completes the proof.
\end{proof}

This lemma yields the desired result for $X(q)$.

\begin{theorem}\label{thm:interp}
The $q$-analogue of the hook length formula for an $n\times n$ square Young diagram is, up to a shift, the major index generating function for reduced expressions of the longest element in the hyperoctahedral group:
\[ \sum_{w \in R(w_0)} q^{\maj(w)} = q^{n\binom{n}{2}}\cdot \frac{[n^2]_q!}{[n]_q^n\prod_{i=1}^{n-1} ([i]_q\cdot [2n-i]_q)^i}.\]
\end{theorem}

Theorem \ref{thm:interp}, along with Corollary \ref{cor: sieving2}, completes the proof of our main result, Theorem \ref{theorem:main}. Because this result can be stated purely in terms of the set $R(w_0)$ and a natural statistic on this set, it would be interesting to obtain a self-contained proof, i.e., one that does not appeal to Haiman's or Rhoades' work. Why must a result about cyclic rotation of words rely on promotion of Young tableaux?

\end{document}